\documentclass[amstex,12pt,russian,amssymb]{article}

\usepackage{mathtext}
\usepackage[cp1251]{inputenc}
\usepackage[T2A]{fontenc}
\usepackage[russian]{babel}
\usepackage[dvips]{graphicx}
\usepackage{amsmath}
\usepackage{amssymb}
\usepackage{amsxtra}
\usepackage{latexsym}
\usepackage{ifthen}

\begin{document}

\newcounter{lemma}
\newcommand{\lemma}{\par \refstepcounter{lemma}%
{\bf Лемма \arabic{lemma}.}}

\newcounter{corollary}
\newcommand{\corollary}{\par \refstepcounter{corollary}%
{\bf Следствие \arabic{corollary}.}}

\newcounter{remark}
\newcommand{\remark}{\par \refstepcounter{remark}%
{\bf Замечание \arabic{remark}.}}

\newcounter{theorem}
\newcommand{\theorem}{\par \refstepcounter{theorem}%
{\bf Теорема \arabic{theorem}.}}

\newcounter{proposition}
\newcommand{\proposition}{\par \refstepcounter{proposition}%
{\bf Предложение \arabic{proposition}.}}

\renewcommand{\refname}{\centerline{\bf Список литературы}}

\newcommand{\proof}{{\it Доказательство.\,\,}}

\noindent УДК 517.5

\medskip
{\bf Е.А.~Севостьянов} (Житомирский государственный университет им.\
И.~Франко)

\medskip
{\bf Є.О.~Севостьянов} (Житомирський державний університет ім.\
І.~Франко)

\medskip
{\bf E.A.~Sevost'yanov} (Zhitomir Ivan Franko State University)

\medskip
{\bf О равностепенной непрерывности отображений с ветвлением в
замыкании области}

\medskip
{\bf Про одностайну неперервність відображень з розгалуженням в
замиканні області}

\medskip
{\bf On equicontinuity of mappings with branching in a closure of a
domain}

\medskip
В настоящей работе изучаются вопрос о локальном поведении
отображений $f:D\rightarrow \overline{{\Bbb R}^n},$ $n\ge 2,$ в
замыкании области $D.$ При определённых условиях на измеримую
функцию $Q(x),$ $Q:D\rightarrow [0, \infty],$ и границы областей $D$
и $D^{\,\prime}=f(D),$ показано, что семейство открытых дискретных
отоб\-ра\-же\-ний $f:D\rightarrow \overline{{\Bbb R}^n},$ имеющих
характеристику квазиконформности $Q(x),$ равностепенно непрерывно в
$\overline{D}.$

\medskip
В даній роботі вивчається питання про локальну поведінку відображень
$f:D\rightarrow \overline{{\Bbb R}^n},$ $n\ge 2,$ в замиканні
області $D.$ За певних умов на вимірну функцію $Q(x),$
$Q:D\rightarrow [0, \infty],$ і межі областей $D$ і
$D^{\,\prime}=f(D),$ показано, що сім'я відкритих дискретних
ві\-доб\-ра\-жень $f:D\rightarrow \overline{{\Bbb R}^n},$ які мають
характеристику квазіконформності $Q(x),$ одностайно неперервна в
$\overline{D}.$

\medskip
In the present paper, questions about a local behavior of mappings
$f:D\rightarrow \overline{{\Bbb R}^n},$ $n\ge 2,$ in $\overline{D}$
are studied. Under some conditions on a measurable function $Q(x),$
$Q:D\rightarrow [0, \infty],$ and boundaries of $D$ and
$D^{\,\prime}=f(D),$ we show that a family of open discrete
map\-ping $f:D\rightarrow \overline{{\Bbb R}^n},$ with
characteristic of quasiconformality $Q(x),$ is equicontinuous in
$\overline{D}.$

\newpage

{\bf 1. Введение.} Всюду далее $D$ -- область в ${\Bbb R}^n,$ $n\ge
2,$ $m$ -- мера Лебега ${\Bbb R}^n,$ запись $f:D\rightarrow {\Bbb
R}^n$ предполагает, что отображение $f,$ заданное в области $D,$
непрерывно. Запись ${\rm dist}(A, B)$ означает евклидово расстояние
между множествами $A$ и $B\subset {\Bbb R}^n.$ Другие определения и
обозначения, встречающиеся в тексте, но не приведённые ниже, могут
быть найдены в работе \cite{Sev$_1$} и монографии \cite{MRSY}.
Граница $\partial D,$ замыкание $\overline{D}$ области $D\subset
{\Bbb R}^n$ (либо области $D\subset \overline{{\Bbb R}^n}$), а также
наличие предела для отображения $f:D\rightarrow {\Bbb R}^n$ (либо
$f:D\rightarrow \overline{{\Bbb R}^n}$) в дальнейшем будут
пониматься в смысле пространства $\overline{{\Bbb R}^n}$
относительно {\it хордальной метрики $h,$} определённой
соотношениями
$$h(x,\infty)=\frac{1}{\sqrt{1+{|x|}^2}}\,, \quad
h(x,y)=\frac{|x-y|}{\sqrt{1+{|x|}^2} \sqrt{1+{|y|}^2}}\,, \quad x\ne
\infty\ne y\,.$$
Хордальным расстоянием $h(A, B)$ между множествами $A, B\subset
\overline{{\Bbb R}^n}$ и хордальным диаметром множества $C\subset
\overline{{\Bbb R}^n}$ называются величины
$$h(A, B):=\inf\limits_{x\in A, y\in B}h(x, y)\,,\quad h(C):=\sup\limits_{x, y\in C}h(x, y)\,,$$
соответственно.

\medskip
В сравнительно недавней работе \cite{Sev$_2$} нами было установлено
свойство равностепенной непрерывности одного семейства
пространственных отображений с неограниченной характеристикой в
предположении, что все отображения рассматриваемого класса являются
гомеоморфизмами. Здесь речь идёт о равностепенной непрерывности в
замыкании $\overline{D}$ области $D\subset {\Bbb R}^n,$ $n\geqslant
2,$ а не только во внутренних точках. По этому поводу см. также
классический результат Някки и Палка для квазиконформных отображений
(см. \cite{NP}). В настоящей заметке будет показано, что
предположение гомеоморфности отображений в \cite{Sev$_2$} может быть
ослаблено до условий открытости и дискретности, при этом, здесь
необходимо требовать условие замкнутости, эквивалентное свойству
сохранения границы (см. \cite[теорема~3.3]{Vu}; см. также
\cite{Zel}), а также и ещё одно условие, заключающееся в следующем:
существует континуум $K\subset D^{\,\prime}=f(D),$ такой что $
h(f^{-1}(K),
\partial D)\geqslant \delta>0$ для некоторого $\delta>0$ и всех
отображений $f$ из рассматриваемого семейства.

\medskip
Приведём теперь некоторые вспомогательные сведения, включая
формулировку основных результатов. Пусть $E,$ $F\subset
\overline{{\Bbb R}^n}$ -- произвольные множества. Обозначим через
$\Gamma(E,F,D)$ семейство всех кривых
$\gamma:[a,b]\rightarrow\overline{{\Bbb R}^n},$ которые соединяют
$E$ и $F$ в $D,$ т.е. $\gamma(a)\in E,$ $\gamma(b)\in F$ и
$\gamma(t)\in D$ при $t\in (a, b).$ Здесь и далее
\begin{equation}\label{eq49***}
A(x_0, r_1,r_2): =\left\{ x\,\in\,{\Bbb R}^n:
r_1<|x-x_0|<r_2\right\}\,.
\end{equation}
Введём в рассмотрение следующее понятие, см. \cite[разд.~7.6
гл.~7]{MRSY}.
Пусть $p\geqslant 1$ и $Q:{\Bbb R}^n\rightarrow [0, \infty]$ --
измеримая по Лебегу функция, $Q(x)\equiv 0$ при всех $x\not\in D.$
Говорят, что отображение $f:D\rightarrow \overline{{\Bbb R}^n}$ есть
{\it кольцевое $Q$-отоб\-ра\-же\-ние в точке $x_0\in \overline{D}$
относительно $p$-модуля,} $x_0\ne \infty,$ если для некоторого
$r_0=r(x_0)$ и произвольных сферического кольца (\ref{eq49***}) и
любых континуумов $E_1\subset \overline{B(x_0, r_1)}\cap D,$
$E_2\subset \left(\overline{{\Bbb R}^n}\setminus B(x_0,
r_2)\right)\cap D,$ отображение $f$ удовлетворяет соотношению
\begin{equation}\label{eq3*!!}
 M_p\left(f\left(\Gamma\left(E_1,\,E_2,\,D\right)\right)\right)\ \le
\int\limits_{A} Q(x)\cdot \eta^p(|x-x_0|)\ dm(x) \end{equation}
для каждой измеримой функции $\eta : (r_1,r_2)\rightarrow [0,\infty
],$ такой что
\begin{equation}\label{eq28*}
\int\limits_{r_1}^{r_2}\eta(r)\ dr\ \ge\ 1\,.
\end{equation}
В точке $x_0=\infty$ данное определение может быть переформулировано
при помощи инверсии: $\varphi(x)=\frac{x}{|x|^2},$ $\infty\mapsto
0.$

Положим $M(\Gamma):=M_n(\Gamma).$ Отметим, что при $p=n$ и $Q(x)\le
K^{\,\prime}$ соотношение (\ref{eq3*!!}) влечёт условие
$M(f(\Gamma))\le K^{\,\prime}\cdot M(\Gamma)$  для семейств кривых
$\Gamma,$ соединяющих сферы $S(x_0, r_1)$ и $S(x_0, r_2),$ как
только $x_0\in D$ и $0<r_1<r_2<{\rm dist}(x_0, \partial D);$ правда
этого, вообще говоря, нельзя сказать относительно любого семейства
$\Gamma$ кривых $\gamma$ в $D.$ Отметим также, что произвольное
отображение $f$ с ограниченным искажением удовлетворяет соотношениям
вида (\ref{eq3*!!})--(\ref{eq28*}) с $Q,$ равным некоторой
постоянной. (Этот результат установлен Е.А. Полецким \cite{Pol}).

\medskip
Компактное множество $G\subset {\Bbb R}^n$ условимся называть {\it
множеством нулевой ёмкости,} пишем ${\rm cap}\,G =0,$ если
существует ограниченное открытое множество $A\subset {\Bbb R}^n,$
такое что $M(\Gamma(G, \partial A, A))=0,$ см., напр.,
\cite[разд.~2, гл.~III и предложение~10.2, гл.~II]{Ri}. Будем
говорить, что произвольное множество $G\subset {\Bbb R}^n$ имеет
ёмкость нуль, если произвольное его компактное подмножество $G_0$
имеет нулевую ёмкость. Условимся также считать, что произвольное
множество $G\subset \overline{{\Bbb R}^n}$ имеет ёмкость нуль, если
${\rm cap}(G\setminus\{\infty\})=0.$

\medskip Область $D$ называется {\it локально связной в точ\-ке}
$x_0\in\partial D,$ если для любой окрестности $U$ точки $x_0$
найдется окрестность $V\subset U$ точки $x_0$ такая, что $V\cap D$
связно, см. \cite[c.~232]{Ku}. Будем говорить, что граница $\partial
D$ области $D$ {\it сильно достижима в точке $x_0\in
\partial D,$} если для любой окрестности $U$ точки $x_0$ найдется
компакт $E\subset D,$ окрестность $V\subset U$ точки $x_0$ и число
$\delta
>0$ такие, что $M(\Gamma(E,F, D))\ge \delta$ для любого континуума  $F$ в $D,$
пересекающего $\partial U$ и $\partial V,$  см., напр.,
\cite[разд.~3.8]{MRSY}.

\medskip
Следуя \cite[раздел 7.22]{He} будем говорить, что борелева функция
$\rho\colon  X\rightarrow [0, \infty]$ является {\it верхним
градиентом} функции $u\colon X\rightarrow {\Bbb R},$ если для всех
спрямляемых кривых $\gamma,$ соединяющих точки $x$ и $y\in X,$
выполняется неравенство $|u(x)-u(y)|\leqslant
\int\limits_{\gamma}\rho\,|dx|,$ где, как обычно,
$\int\limits_{\gamma}\rho\,|dx|$ обозначает линейный интеграл от
функции $\rho$ по кривой $\gamma.$ Будем также говорить, что в
указанном пространстве $X$ выполняется {\it $(1; p)$-неравенство
Пуанкаре,} если найдётся постоянная $C\geqslant 1$ такая, что для
каждого шара $B\subset X,$ произвольной ограниченной непрерывной
функции $u\colon X\rightarrow {\Bbb R}$ и любого её верхнего
градиента $\rho$ выполняется следующее неравенство:
$$\frac{1}{\mu(B)}\int\limits_{B}|u-u_B|d\mu(x)\leqslant C\cdot({\rm diam\,}B)\left(\frac{1}{\mu(B)}
\int\limits_{B}\rho^p d\mu(x)\right)^{1/p}\,,$$
где $u_B:=\frac{1}{\mu(B)}\int\limits_Bu(x)d\mu(x).$ Метрическое
пространство $(X, d, \mu)$ назовём {\it $n$-регулярным по Альфорсу,}
если при каждом $x_0\in X,$ некоторой постоянной $C\geqslant 1$ и
всех $R<{\rm diam}\,X$
$$\frac{1}{C}R^{n}\leqslant \mu(B(x_0, R))\leqslant CR^{n}\,.$$

\begin{remark}\label{rem1}
Одним из примеров $n$-регулярного по Альфорсу пространства
относительно евклидовой метрики и меры Лебега в ${\Bbb R}^n,$ в
котором выполнено $(1; p)$-неравенство Пуанкаре, является единичный
шар ${\Bbb B}^n,$ см. \cite[предложение~2.1]{Sev$_6$}.
\end{remark}

\medskip Согласно \cite[разд.~6.1, гл.~6]{MRSY}, будем говорить, что
функция ${\varphi}:D\rightarrow{\Bbb R}$ имеет {\it конечное среднее
колебание} в точке $x_0\in D$, пишем $\varphi\in FMO(x_0),$ если
$$\overline{\lim\limits_{\varepsilon\rightarrow 0}}\quad\frac{1}{\Omega_n\varepsilon^n}
\int\limits_{B(x_0, \varepsilon)}
|{\varphi}(x)-\overline{{\varphi}}_{\varepsilon}|\,dm(x)<\infty\,,
$$
где $\Omega_n$ -- объём единичного шара ${\Bbb B}^n$ в ${\Bbb R}^n$
и
$\overline{{\varphi}}_{\varepsilon}\,=\,\frac{1}{\Omega_n\varepsilon^n}\int\limits_{B(
x_0,\,\varepsilon)} {\varphi}(x)\,dm(x).$

\medskip
Для $p\geqslant 1,$ фиксированных областей $D\subset {\Bbb R}^n$ и
$D^{\,\prime}$ в $\overline{{\Bbb R}^n},$ $n\geqslant 2,$ континуума
$K\subset D^{\,\prime}$ и числа $\delta>0$ обозначим через
$\frak{F}_{Q, \delta, K, p}(D, D^{\,\prime})$ семейство всех
открытых дискретных кольцевых $Q$-отображений $f:D\rightarrow
D^{\,\prime}$ относительно $p$-модуля, таких что
$f(D)=D^{\,\prime},$ и для которых $h(f^{-1}(K),
\partial D)\geqslant \delta>0.$

\medskip
Пусть $Q:D\rightarrow [0,\infty]$ -- измеримая по Лебегу функция,
тогда $q_{x_0}(r)$ обозначает среднее интегральное значение $Q(x)$
над сферой $|x-x_0|=r,$ $
q_{x_0}(r):=\frac{1}{\omega_{n-1}r^{n-1}}\int\limits_{|x-x_0|=r}Q(x)\,dS,$
где $dS$ -- элемент площади поверхности $S.$ Полагаем
\begin{equation*}\label{eq8} q^{\,\prime}_{b}(r):=\frac{1}{\omega_{n-1}r^{n-1}}\int\limits_{|x-b|=r}Q^{\,\prime}(x)\,dS\,,
\quad Q^{\,\prime}(x)=\left\{
\begin{array}{rr}
Q(x), &   Q(x)\ge 1\,, \\
1,  &  Q(x)<1\,.
\end{array}
\right. \end{equation*}
Имеют место следующие утверждения.

\medskip
\begin{theorem}\label{th3} {\sl\, Предположим, область $D$ локально связна в каждой точке
$b\in\partial D,$ $C(f, \partial D)\subset D^{\,\prime}$ для каждого
$f\in\frak{F}_{Q, \delta, K, n}(D, D^{\,\prime}),$ и что область
$D^{\,\prime}$ имеет сильно достижимую границу. Если функция $Q$
имеет конечное среднее колебание в $\overline{D},$ либо в каждой
точке $x_0\in \overline{D}$
$$\int\limits_{0}^{\delta(x_0)}\frac{dt}{tq_{x_0}^{\,\prime\,\frac{1}{n-1}}(t)}=\infty\,,$$
то каждое из отображений $f\in\frak{F}_{Q, \delta, K, n}(D,
D^{\,\prime})$ имеет непрерывное продолжение в $\overline{D}.$ Если,
кроме того, ${\rm cap}\,({\Bbb R}^n\setminus D^{\,\prime})>0,$ то
семейство $\frak{F}_{Q, \delta, K, n}(D, D^{\,\prime}),$ состоящее
из всех, таким образом, продолженных отображений $\overline{f}:
\overline{D}\rightarrow \overline{D^{\,\prime}},$ является
равностепенно непрерывным в $\overline{D}.$
  }
\end{theorem}

\medskip
\begin{theorem}\label{th4} {\sl\, Предположим, $n-1<p\leqslant n,$ область $D$ локально связна в каждой точке
$b\in\partial D,$ $C(f, \partial D)\subset D^{\,\prime}$ для каждого
$f\in\frak{F}_{Q, \delta, K, p}(D, D^{\,\prime}),$ и что область
$D^{\,\prime}\subset {\Bbb R}^n$ ограничена и является
$n$-регулярным по Альфорсу пространством относительно евклидовой
метрики и меры Лебега в ${\Bbb R}^n,$ в котором выполнено $(1;
p)$-неравенство Пуанкаре. Если функция $Q$ имеет конечное среднее
колебание в $\overline{D},$ либо в каждой точке $x_0\in
\overline{D}$
$$\int\limits_{0}^{\delta(x_0)}\frac{dt}{t^{\frac{n-1}{p-1}}q_{x_0}^{\,\prime\,\frac{1}{p-1}}(t)}=\infty\,,$$
то каждое из отображений $f\in\frak{F}_{Q, \delta, K, p}(D,
D^{\,\prime})$ имеет непрерывное продолжение в $\overline{D}$ и
семейство $\frak{F}_{Q, \delta, K, p}(D, D^{\,\prime}),$ состоящее
из всех, таким образом, продолженных отображений $\overline{f}:
\overline{D}\rightarrow \overline{D^{\,\prime}},$ является
равностепенно непрерывным в $\overline{D}.$ }
\end{theorem}

\medskip
{\bf 2. Вспомогательные сведения.} Справедливо следующее утверждение
(см.~\cite[предложение~4.7]{AS}).

\medskip
 \begin{proposition}\label{pr2}
{\sl Пусть $X$~--- $n$-регулярное по Альфорсу метрическое
пространство с мерой\/{\em,} в котором выполняется
$(1;p)$-неравенство Пуанкаре. Тогда для произвольных континуумов $E$
и $F,$ содержащихся в шаре $B(x_0, R),$ и некоторой постоянной $C>0$
выполняется неравенство
$$M_p(\Gamma(E, F, X))\geqslant \frac{1}{C}\cdot\frac{\min\{{\rm diam}\,E, {\rm diam}\,F\}}{R^{1+p-n}}\,.$$ }
\end{proposition}

Здесь величина $M_p(\Gamma(E, F, X))$ может быть определена по
аналогии с евклидовым случаем (см. \cite[стр.~610]{AS}). Аналог
следующего утверждения при $p=n$ установлен в
\cite[лемма~1]{Sev$_4$}. Однако, ниже мы доказываем это утверждение
при несколько иных условиях, связанных с регулярностью по Альфорсу и
неравенством типа Пуанкаре.

\medskip
\begin{lemma}\label{pr1}
{\sl\, Пусть $n-1\leqslant p\leqslant n,$ $f:D\rightarrow {\Bbb
R}^n$ -- открытое дискретное кольцевое $Q$-отоб\-ра\-же\-ние в точке
$b\in
\partial D$ относительно $p$-модуля, $f(D)=D^{\,\prime},$ область $D$ локально связна в
точке $b,$ $C(f,
\partial D)\subset \partial D^{\,\prime},$ область $D^{\,\prime}$
ограничена и является $n$-регулярным по Альфорсу пространством
относительно евклидовой метрики и меры Лебега в ${\Bbb R}^n,$ в
котором выполнено $(1; p )$-неравенство Пуанкаре. Предположим, что
найдётся $\varepsilon_0>0$ и некоторая неотрицательная измеримая по
Лебегу функция $\psi(t),$ $\psi:(0, \varepsilon_0)\rightarrow
[0,\infty],$ такая что для всех $\varepsilon\in(0, \varepsilon_0)$
\begin{equation}\label{eq7***} I(\varepsilon,
\varepsilon_0):=\int\limits_{\varepsilon}^{\varepsilon_0}\psi(t)dt <
\infty\,,\quad I(\varepsilon, \varepsilon_0)\rightarrow
\infty\quad\text{при}\quad\varepsilon\rightarrow 0
\end{equation}
и при $\varepsilon\rightarrow 0$
\begin{equation}\label{eq5***}
\int\limits_{A(\varepsilon, \varepsilon_0, b)}
Q(x)\cdot\psi^{\,p}(|x-b|)
 \ dm(x) =o(I^{p}(\varepsilon, \varepsilon_0))\,,
\end{equation}
где $A:=A(\varepsilon, \varepsilon_0, b)$ определено в
(\ref{eq49***}). Тогда $C(f, b)=\{y\}.$}
\end{lemma}

\medskip
\begin{proof}
Предположим противное. Тогда найдутся, по крайней мере, две
последовательности $x_i,$ $x_i^{\,\prime}\in D,$ $i=1,2,\ldots,$
такие, что $x_i\rightarrow b,$ $x_i^{\,\prime}\rightarrow b$ при
$i\rightarrow \infty,$ $f(x_i)\rightarrow y,$
$f(x_i^{\,\prime})\rightarrow y^{\,\prime}$ при $i\rightarrow
\infty$ и $y^{\,\prime}\ne y.$ Отметим, что $y$ и $y^{\,\prime}\in
\partial D^{\,\prime},$ поскольку по условию $C(f,
\partial D)\subset \partial D^{\,\prime}.$ Так как область $D^{\,\prime}$
ограничена и является $n$-регулярным по Альфорсу пространством
относительно евклидовой метрики и меры Лебега в ${\Bbb R}^n,$ в
котором выполнено $(1; p )$-неравенство Пуанкаре, то ввиду
предложения \ref{pr2}
\begin{equation}\label{eq1}
M_p(\Gamma(C_0^{\,\prime}, F, D^{\,\prime}))\geqslant
\frac{1}{C}\cdot\frac{\delta}{R^{1+p-n}}
>0\,,
\end{equation}
где $R>0$ таково, что $D^{\,\prime}\subset B(x_0, R),$ $x_0$ --
некоторая фиксированная точка области $D^{\,\prime},$ а $F$ и
$C_0^{\,\prime}$ -- произвольные фиксированные континуумы в
$D^{\,\prime},$ диаметры которых не меньше $\delta.$ Поскольку
область $D$ локально связна в точке $b,$ можно соединить точки $x_i$
и $x_i^{\,\prime}$ кривой $\gamma_i,$ лежащей в $V\cap D.$ Можно
также считать, что $|\gamma_i|\in \overline{B(b, 2^{\,-i})}\cap D,$
где $|\gamma_i|:=\{x\in {\Bbb R}^n: \exists\,t: \gamma_i(t)=x\}.$
Поскольку $f(x_i)\rightarrow y,$ $f(x_i^{\,\prime})\rightarrow
y^{\,\prime}$ при $i\rightarrow \infty$ и $y^{\,\prime}\ne y,$
существует $\delta>0$ такое, что ${\rm diam\,}|f(\gamma_i)|\geqslant
\delta$ при всех $i\in {\Bbb N}$ (здесь $|f(\gamma_i)|:=\{x\in {\Bbb
R}^n: \exists\,t: f(\gamma_i(t))=x\}$). Пусть $C_0^{\,\prime}$ --
произвольный континуум в $D^{\,\prime},$ также имеющий диаметр не
меньший $\delta,$ тогда ввиду (\ref{eq1})
\begin{equation}\label{eq2}
M_p(\Gamma(C_0^{\,\prime}, |f(\gamma_i)|, D^{\,\prime}))\geqslant
\frac{1}{C}\cdot\frac{\delta}{R^{1+p-n}}
\end{equation}
при всех $i\in {\Bbb N}.$

Обозначим через $\Gamma_i$ семейство всех полуоткрытых кривых
$\beta_i:[a, b)\rightarrow {\Bbb R}^n$ таких, что $\beta_i(a)\in
|f(\gamma_i)|,$ $\beta_i(t)\in D^{\,\prime}$ при всех $t\in [a, b)$
и, кроме того, $\lim\limits_{t\rightarrow b-0}\beta_i(t):=B_i\in
C_0^{\,\prime}.$ Очевидно, что
\begin{equation}\label{eq4}
M_p(\Gamma_i)=M_p\left(\Gamma\left(C_0^{\,\prime}, |f(\gamma_i)|,
D^{\,\prime}\right)\right)\,.
\end{equation}
При каждом фиксированном $i\in {\Bbb N},$ $i\ge i_0,$ рассмотрим
семейство $\Gamma_i^{\,\prime}$ максимальных поднятий
$\alpha_i(t):[a, c)\rightarrow D$ семейства $\Gamma_i$ с началом во
множестве $|\gamma_i|.$ Такое семейство существует и определено
корректно ввиду \cite[следствие~3.3, гл. II]{Ri}. Заметим, прежде
всего, что никакая кривая $\alpha_i(t)\in \Gamma_i^{\,\prime},$
$\gamma_i:[a, c)\rightarrow D,$ не может стремиться к границе
области $D$ при $t\rightarrow c-0$ ввиду условия $C(f,
\partial D)\subset \partial D^{\,\prime}.$ Тогда $C(\alpha_i(t), c)\subset
D.$ Предположим теперь, что кривая $\alpha_i(t)$ не имеет предела
при $t\rightarrow c-0.$ Тогда предельное множество $C(\alpha_i(t),
c)$ есть континуум в $D.$ В силу непрерывности отображения $f,$
получаем, что $f\equiv const$ на $C(\alpha_i(t), c),$ что
противоречит предположению о дискретности $f.$

Следовательно, $\exists \lim\limits_{t\rightarrow
c-0}\alpha_i(t)=A_i\in D.$ Отметим, что, в этом случае, по
определению максимального поднятия, $c=b.$ Тогда, с одной стороны,
$\lim\limits_{t\rightarrow b-0}\alpha_i(t):=A_i,$ а с другой, в силу
непрерывности отображения $f$ в $D,$
$$f(A_i)=\lim\limits_{t\rightarrow b-0}f(\alpha_i(t))=\lim\limits_{t\rightarrow b-0}
\beta_i(t)=B_i\in C_0^{\,\prime}\,.$$ Отсюда, по определению $C_0,$
следует, что $A_i\in C_0.$ Погрузим компакт $C_0$ в некоторый
континуум $C_1,$ всё ещё полностью лежащий в области $D,$ что
возможно ввиду \cite[лемма~1]{Smol}. За счёт уменьшения
$\varepsilon_0>0,$ можно снова считать, что $C_1\cap\overline{B(b,
\varepsilon_0)}=\varnothing.$
Заметим, что функция
$$\eta(t)=\left\{
\begin{array}{rr}
\psi(t)/I(2^{-i}, \varepsilon_0), &   t\in (2^{-i},
\varepsilon_0),\\
0,  &  t\in {\Bbb R}\setminus (2^{-i}, \varepsilon_0)\,,
\end{array}
\right. $$ где $I(\varepsilon,
\varepsilon_0)=\int\limits_{\varepsilon}^{\varepsilon_0}\psi(t)dt,$
удовлетворяет условию нормировки вида (\ref{eq28*}) при
$r_1:=2^{-i},$ $r_2:=\varepsilon_0,$ поэтому, в силу определения
кольцевого $Q$-отоб\-ра\-же\-ния в граничной точке относительно
$p$-модуля, а также ввиду условий (\ref{eq7***})--(\ref{eq5***}),
\begin{equation}\label{eq11*}
M_p\left(f\left(\Gamma_i^{\,\prime}\right)\right)\le \Delta(i)\,,
\end{equation}
где $\Delta(i)\rightarrow 0$ при $i\rightarrow \infty.$ Однако,
$\Gamma_i=f(\Gamma_i^{\,\prime}),$ поэтому из (\ref{eq11*}) получим,
что при $i\rightarrow \infty$
\begin{equation}\label{eq3}
M_p(\Gamma_i)= M_p\left(f(\Gamma_i^{\,\prime})\right)\le
\Delta(i)\rightarrow 0\,.
\end{equation}
Однако, соотношение (\ref{eq3}) вместе с равенством (\ref{eq4})
противоречат неравенству (\ref{eq2}), что и доказывает лемму.
\end{proof}$\Box$

\medskip
{\bf 3. Формулировка и доказательство основных лемм.} Основным
инструментом на пути доказательства теорем \ref{th3} и \ref{th4}
являются следующие два утверждения.

\medskip
\begin{lemma}\label{lem1}{\sl\,
Предположим, область $D$ локально связна в каждой точке
$b\in\partial D,$ $C(f, \partial D)\subset D^{\,\prime}$ для каждого
$f\in \frak{F}_{Q, \delta, K, n}(D, D^{\,\prime}),$ а область
$D^{\,\prime}$ имеет сильно достижимую границу. Предположим также,
что для каждой точки $x_0\in \overline{D}$ найдётся
$\varepsilon_0=\varepsilon_0(x_0)<{\rm dist\,}(x_0,\partial D)$ и
измеримая по Лебегу функция $\psi(t):(0, \varepsilon_0)\rightarrow
[0,\infty]$ со следующим свойством: для любого $\varepsilon\in(0,
\varepsilon_0)$
\begin{equation} \label{eq3.7.1}
I(\varepsilon,
\varepsilon_0):=\int\limits_{\varepsilon}^{\varepsilon_0}\psi(t)dt <
\infty\,,\quad I(\varepsilon, \varepsilon_0)\rightarrow
\infty\quad\text{при}\quad\varepsilon\rightarrow 0
\end{equation}
и, кроме того, при  $\varepsilon\rightarrow 0$
\begin{equation} \label{eq3.7.2}
\int\limits_{A(x_0, \varepsilon, \varepsilon_0)}
Q(x)\cdot\psi^{\,n}(|x-x_0|)\,dm(x) = o(I^n(\varepsilon,
\varepsilon_0))\,,\end{equation}
где, как обычно, сферическое кольцо $A(x_0, \varepsilon,
\varepsilon_0)$ определено как в (\ref{eq49***}).
Тогда каждое отображение $\frak{F}_{Q, \delta, K, n}(D,
D^{\,\prime})$ продолжается до непрерывного отображения
$\overline{f}: \overline{D}\rightarrow \overline{D^{\,\prime}}.$
Если, кроме того, ${\rm cap}\,({\Bbb R}^n\setminus D^{\,\prime})>0,$
то семейство $\overline{\frak{F}_{Q, \delta, K, n}(D,
D^{\,\prime})}$ состоящее из всех, таким образом, продолженных
отображений $\overline{f}: \overline{D}\rightarrow
\overline{D^{\,\prime}},$ является равностепенно непрерывным в
$\overline{D}.$}
\end{lemma}

\medskip
\begin{proof}
Равностепенная непрерывность внутри области $D$ вытекает из
\cite[лемма~5.2]{Sev$_3$}, а возможность продолжения каждого
элемента $f$ семейства отображений $\frak{F}_{Q, \delta, K, n}(D,
D^{\,\prime})$ до непрерывного отображения в замыкании $D$ ~--- из
\cite[лемма~1]{Sev$_4$}.

Осталось показать, что семейство $\frak{F}_{Q, \delta, K, n}(D,
D^{\,\prime})$ (обозначения не меняем) равностепенно не\-прерывно в
точках $\partial D.$ Предположим противное, тогда найдётся $x_0\in
\partial D$ и число $a>0$ такое, что для каждого
$m=1,2,\ldots$ существуют точка $x_m\in \overline{D}$ и элемент
${f}_m$ семейства $\frak{F}_{Q, \delta, K, n}(D, D^{\,\prime})$
такие, что $|x_0-x_m|< 1/m$ и
\begin{equation}\label{eq6***}
h\left(f_m(x_m), f_m(x_0)\right)\ge a\,.
\end{equation}
%
Можно считать, что $x_0\ne \infty.$ В виду возможности непрерывного
продолжения каждого $f_m$ на границу $D,$ мы можем считать, что
$x_m\in D.$

В силу локальной связности области $D$ в точке $x_0$ найдётся
последовательность окрестностей $V_m$ точки $x_0$ с ${\rm
diam}\,V_m\rightarrow 0$ при $m\rightarrow\infty,$ такие что
множества $D\cap V_m $ являются областями и $x_m\in D\cap V_m.$ Т.к.
граничные точки области, локально связной на границе являются
достижимыми из $D$ некоторым локально спрямляемым путём, см.
\cite[предложение~13.2]{MRSY}, мы можем соединить точки $x_m$ и
$x_0$ непрерывной кривой $\gamma_m(t):[0,1]\rightarrow {\Bbb R}^n$
такой, что $\gamma_m(0)=x_0,$ $\gamma_m(1)=x_m$ и $\gamma_m(t)\in
V_m$ при $t\in (0,1).$ Обозначим через $C_m$ образ кривой
$\gamma_m(t)$ при отображении $f_m.$ Из соотношения (\ref{eq6***})
вытекает, что
\begin{equation}\label{eq5.1}
h(C_m)\geqslant a\qquad\forall\, m\in {\Bbb N}\,,
\end{equation}
где $h$ обозначает хордальный диаметр множества. Поскольку $f_m$
сохраняет границу, точка $f_m(x_0)$ принадлежит $\partial
D^{\,\prime}.$ Поскольку $\overline{{\Bbb R}^n}$ компактно, не
ограничивая общности можно считать, что последовательность
$f_m(x_0)$ сходится к некоторой точке $y_0\in
\partial D^{\,\prime}$ при $m\rightarrow\infty.$
Заметим, что если $E$ -- компакт из определения сильно достижимой
границы области $D^{\,\prime},$ то вместо него может быть взят
произвольный континуум из $D^{\,\prime}$ (см.
\cite[лемма~4.1]{Sev$_2$}). Таким образом, из определения сильно
достижимой границы в точке $y_0,$ с учётом условия (\ref{eq5.1})
найдется $b>0$ такое что
\begin{equation}\label{eq13}
M(\Gamma(K, C_m, D^{\,\prime}))\geqslant b\qquad\forall\, m\in {\Bbb
N}\,.
\end{equation}
С другой стороны, рассмотрим теперь семейство $\Gamma_m^1,$
состоящее из всех максимальных поднятий $\alpha:[0, c)\rightarrow D$
семейства $\Gamma_m:=\Gamma(K, C_m, D^{\,\prime})$ при отображении
$f_m$ с началом в $|\gamma_m|=\{x\in D: \exists\, t:
\gamma_m(t)=x\}.$ Поскольку все отображения $f_m$ являются открытыми
и дискретными и, кроме того, $C(f_m,
\partial D)\subset \partial D^{\,\prime}$ при каждом $m\in {\Bbb N},$
указанное семейство максимальных поднятий существует и
$\Gamma_m^1\subset \Gamma(|\gamma_m|, f_m^{\,-1}(K), D).$ Согласно
сказанному,
 Поскольку
при каждом фиксированном $m\in {\Bbb N}$ множество $|\gamma_m|$
принадлежит окрестности $V_m$ точки $x_0,$ где ${\rm
diam\,}V_m\rightarrow 0$ при $m\rightarrow \infty,$ для
последовательности $\varepsilon_k=\frac{1}{2^k}$ найдётся
подпоследовательность номеров $m_k,$ $k=1,2,\ldots ,$ таких что
$\gamma_{m_k}\subset B(x_0, \frac{1}{2^k}).$ Заметим, что ввиду
компактности пространства $\overline{{\Bbb R}^n}$ при каждом
фиксированном $\delta>0$ множество
$$C_{\delta}:=\{x\in D: h(x, \partial D)\geqslant \delta\}$$
является компактом в $D$ и $f_{m_k}^{-1}(K)\subset C_{\delta}$ при
некотором $\delta>0$ и всех натуральных $k.$ Ввиду
\cite[лемма~1]{Smol} множество $C_{\delta}$ можно вложить в
континуум $E_{\delta},$ лежащий в области $D,$ при этом, можно
считать, что ${\rm dist}\,(x_0, E_{\delta})\geqslant \varepsilon_0$
за счёт уменьшения $\varepsilon_0,$ если это необходимо. Тогда на
основании (\ref{eq3*!!}) вытекает, что
\begin{equation}\label{eq10}
M(f_{m_k}(\Gamma_{m_k}^1))\leqslant M(f_{m_k}(\Gamma(|\gamma_{m_k}|,
E_{\delta}, D)))\le \int\limits_{A(x_0, \frac{1}{2^k},
\varepsilon_0)} Q(x)\cdot \eta^n(|x-x_0|)\ dm(x)
\end{equation}
для каждой измеримой функции $\eta: (\frac{1}{2^k},
\varepsilon_0)\rightarrow [0,\infty ],$ такой что
$\int\limits_{\frac{1}{2^k}}^{\varepsilon_0}\eta(r)dr \geqslant 1.$
Заметим, что функция
$$\eta(t)=\left\{
\begin{array}{rr}
\psi(t)/I(2^{-k}, \varepsilon_0), &   t\in (2^{-k},
\varepsilon_0),\\
0,  &  t\in {\Bbb R}\setminus (2^{-k}, \varepsilon_0)\,,
\end{array}
\right. $$ где $I(\varepsilon,
\varepsilon_0):=\int\limits_{\varepsilon}^{\varepsilon_0}\psi(t)dt,$
удовлетворяет условию нормировки вида (\ref{eq28*}) при
$r_1:=2^{-k},$ $r_2:=\varepsilon_0,$ поэтому из условий (\ref{eq10})
и (\ref{eq3.7.2}) вытекает, что
\begin{equation}\label{eq11}
M(f_{m_k}(\Gamma_{m_k}^1))\leqslant \alpha(2^{-k})\rightarrow 0
\end{equation}
при $k\rightarrow \infty,$ где $\alpha(\varepsilon)$ -- некоторая
неотрицательная функция, стремящаяся к нулю при
$\varepsilon\rightarrow 0,$ которая существует ввиду условия
(\ref{eq3.7.2}). Заметим, кроме того, что $f(\Gamma_{m_k}^1)>
\Gamma_{m_k}$ и, одновременно,
$f(\Gamma_{m_k}^1)\subset\Gamma_{m_k},$ так что ввиду \cite[теоремы
6.2, 6.4]{Va}
\begin{equation}\label{eq12}
M(f_{m_k}(\Gamma_{m_k}^1))=M(\Gamma_{m_k})\,.
\end{equation}
Однако, соотношения (\ref{eq11}) и (\ref{eq12}) в совокупности
противоречат (\ref{eq13}). Полученное противоречие указывает на то,
что исходное предположение (\ref{eq6***}) было неверным, и, значит,
семейство отображений $\frak{F}_{Q, \delta, K, n}(D, D^{\,\prime})$
равностепенно непрерывно в каждой точке $x_0\in \partial D.$
\end{proof}$\Box$

\medskip
\begin{lemma}\label{lem1A}{\sl\,
Предположим, $n-1<p\leqslant n,$ область $D$ локально связна в
каждой точке $b\in\partial D,$ а область $D^{\,\prime}\subset {\Bbb
R}^n$ ограничена и является $n$-регулярным по Альфорсу пространством
относительно евклидовой метрики и меры Лебега в ${\Bbb R}^n,$ в
котором выполнено $(1; p )$-неравенство Пуанкаре. Предположим также,
что для каждой точки $x_0\in \overline{D}$ найдётся
$\varepsilon_0=\varepsilon_0(x_0)<{\rm dist\,}(x_0,\partial D)$ и
измеримая по Лебегу функция $\psi(t):(0, \varepsilon_0)\rightarrow
[0,\infty]$ со следующим свойством: для любого $\varepsilon\in(0,
\varepsilon_0)$
\begin{equation} \label{eq3.7.1A}
I(\varepsilon,
\varepsilon_0):=\int\limits_{\varepsilon}^{\varepsilon_0}\psi(t)dt <
\infty\,,\quad I(\varepsilon, \varepsilon_0)\rightarrow
\infty\quad\text{при}\quad\varepsilon\rightarrow 0
\end{equation}
и, кроме того, при  $\varepsilon\rightarrow 0$
\begin{equation} \label{eq3.7.2A}
\int\limits_{A(x_0, \varepsilon, \varepsilon_0)}
Q(x)\cdot\psi^{\,p}(|x-x_0|)\,dm(x) = o(I^p(\varepsilon,
\varepsilon_0))\,,\end{equation}
где, как обычно, сферическое кольцо $A(x_0, \varepsilon,
\varepsilon_0)$ определено как в (\ref{eq49***}).
Тогда каждое отображение $\frak{F}_{Q, \delta, K, p}(D,
D^{\,\prime})$ продолжается до непрерывного отображения
$\overline{f}: \overline{D}\rightarrow \overline{D^{\,\prime}}$ и
семейство $\overline{\frak{F}_{Q, \delta, K, p}(D, D^{\,\prime})}$
состоящее из всех, таким образом, продолженных отображений
$\overline{f}: \overline{D}\rightarrow \overline{D^{\,\prime}},$
является равностепенно непрерывным в $\overline{D}.$}
\end{lemma}

\medskip
\begin{proof} Схема доказательства полностью аналогична схеме,
изложенной при доказательстве леммы \ref{lem1}, однако, ради полноты
изложения мы приведём его полностью. Равностепенная непрерывность
внутри области $D$ вытекает из \cite[лемма~5.2]{Sev$_3$} при $p=n$ и
\cite[лемма~3]{Sev$_5$}, а возможность продолжения каждого элемента
$f$ семейства отображений $\frak{F}_{Q, \delta, K, p}(D,
D^{\,\prime})$ до непрерывного отображения в замыкании $D$ ~--- из
леммы \ref{pr1}.

Осталось показать, что семейство $\frak{F}_{Q, \delta, K, p}(D,
D^{\,\prime})$ (обозначения не меняем) равностепенно не\-прерывно в
точках $\partial D.$ Предположим противное, тогда найдётся $x_0\in
\partial D$ и число $a>0$ такое, что для каждого
$m=1,2,\ldots$ существуют точка $x_m\in \overline{D}$ и элемент
${f}_m$ семейства $\frak{F}_{Q, \delta, K, p}(D, D^{\,\prime})$
такие, что $|x_0-x_m|< 1/m$ и
\begin{equation}\label{eq6***A}
h\left(f_m(x_m), f_m(x_0)\right)\ge a\,.
\end{equation}
%
В виду возможности непрерывного продолжения каждого $f_m$ на границу
$D,$ мы можем считать, что $x_m\in D.$

В силу локальной связности области $D$ в точке $x_0$ найдётся
последовательность окрестностей $V_m$ точки $x_0$ с ${\rm
diam}\,V_m\rightarrow 0$ при $m\rightarrow\infty,$ такие что
множества $D\cap V_m $ являются областями и $x_m\in D\cap V_m.$ Т.к.
граничные точки области, локально связной на границе являются
достижимыми из $D$ некоторым локально спрямляемым путём, см.
\cite[предложение~13.2]{MRSY}, мы можем соединить точки $x_m$ и
$x_0$ непрерывной кривой $\gamma_m(t):[0,1]\rightarrow {\Bbb R}^n$
такой, что $\gamma_m(0)=x_0,$ $\gamma_m(1)=x_m$ и $\gamma_m(t)\in
V_m$ при $t\in (0,1).$ Обозначим через $C_m$ образ кривой
$\gamma_m(t)$ при отображении $f_m.$ Из соотношения (\ref{eq6***})
вытекает, что
\begin{equation}\label{eq5.1A}
d(C_m)\geqslant h(C_m)\geqslant a\qquad\forall\, m\in {\Bbb N}\,,
\end{equation}
где $h$ обозначает хордальный диаметр множества, а $d$ -- его
евклидов диаметр. Поскольку область $D^{\,\prime}$ ограничена,
существует шар $B_R=\{x\in D^{\,\prime}: |x-\overline{x_0}|<R\}$ с
центром в некоторой точке $\overline{x_0}\in D^{\,\prime},$
совпадающий с $D^{\,\prime}.$ Тогда из предложения \ref{pr2} с
учётом неравенства (\ref{eq5.1A}) следует, что
\begin{equation}\label{eq13A}
M_p(\Gamma(K, C_m, D^{\,\prime}))\geqslant b\qquad\forall\, m\in
{\Bbb N}\,.
\end{equation}
С другой стороны, рассмотрим теперь семейство $\Gamma_m^1,$
состоящее из всех максимальных поднятий $\alpha:[0, c)\rightarrow D$
семейства $\Gamma_m:=\Gamma(K, C_m, D^{\,\prime})$ при отображении
$f_m$ с началом в $|\gamma_m|=\{x\in D: \exists\, t:
\gamma_m(t)=x\}.$ Поскольку все отображения $f_m$ являются открытыми
и дискретными и, кроме того, $C(f_m,
\partial D)\subset \partial D^{\,\prime}$ при каждом $m\in {\Bbb N},$
указанное семейство максимальных поднятий существует и
$\Gamma_m^1\subset \Gamma(|\gamma_m|, f_m^{\,-1}(K), D)$ (см.
\cite[следствие~3.3, гл. II]{Ri}). Поскольку при каждом
фиксированном $m\in {\Bbb N}$ множество $|\gamma_m|$ принадлежит
окрестности $V_m$ точки $x_0,$ где ${\rm diam\,}V_m\rightarrow 0$
при $m\rightarrow \infty,$ для последовательности
$\varepsilon_k=\frac{1}{2^k}$ найдётся подпоследовательность номеров
$m_k,$ $k=1,2,\ldots ,$ таких что $\gamma_{m_k}\subset B(x_0,
\frac{1}{2^k}).$ Заметим, что ввиду компактности пространства
$\overline{{\Bbb R}^n}$ при каждом фиксированном $\delta>0$
множество
$$C_{\delta}:=\{x\in D: h(x, \partial D)\geqslant \delta\}$$
является компактом в $D$ и $f_{m_k}^{-1}(K)\subset C_{\delta}$ при
некотором $\delta>0$ и всех натуральных $k.$ Ввиду
\cite[лемма~1]{Smol} множество $C_{\delta}$ можно вложить в
континуум $E_{\delta},$ лежащий в области $D,$ при этом, можно
считать, что ${\rm dist}\,(x_0, E_{\delta})\geqslant \varepsilon_0$
за счёт уменьшения $\varepsilon_0,$ если это необходимо. Тогда на
основании (\ref{eq3*!!}) вытекает, что
\begin{equation}\label{eq10A}
M_p(f_{m_k}(\Gamma_{m_k}^1))\leqslant
M_p(f_{m_k}(\Gamma(|\gamma_{m_k}|, E_{\delta}, D)))\le
\int\limits_{A(x_0, \frac{1}{2^k}, \varepsilon_0)} Q(x)\cdot
\eta^p(|x-x_0|)\ dm(x)
\end{equation}
для каждой измеримой функции $\eta: (\frac{1}{2^k},
\varepsilon_0)\rightarrow [0,\infty ],$ такой что
$\int\limits_{\frac{1}{2^k}}^{\varepsilon_0}\eta(r)dr \geqslant 1.$
Заметим, что функция
$$\eta(t)=\left\{
\begin{array}{rr}
\psi(t)/I(2^{-k}, \varepsilon_0), &   t\in (2^{-k},
\varepsilon_0),\\
0,  &  t\in {\Bbb R}\setminus (2^{-k}, \varepsilon_0)\,,
\end{array}
\right. $$ где $I(\varepsilon,
\varepsilon_0):=\int\limits_{\varepsilon}^{\varepsilon_0}\psi(t)dt,$
удовлетворяет условию нормировки вида (\ref{eq28*}) при
$r_1:=2^{-k},$ $r_2:=\varepsilon_0,$ поэтому из условий
(\ref{eq10A}) и (\ref{eq3.7.2A}) вытекает, что
\begin{equation}\label{eq11A}
M_p(f_{m_k}(\Gamma_{m_k}^1))\leqslant \alpha(2^{-k})\rightarrow 0
\end{equation}
при $k\rightarrow \infty,$ где $\alpha(\varepsilon)$ -- некоторая
неотрицательная функция, стремящаяся к нулю при
$\varepsilon\rightarrow 0,$ которая существует ввиду условия
(\ref{eq3.7.2A}). Заметим, кроме того, что $f(\Gamma_{m_k}^1)>
\Gamma_{m_k}$ и, одновременно,
$f(\Gamma_{m_k}^1)\subset\Gamma_{m_k},$ так что ввиду \cite[теоремы
6.2, 6.4]{Va}
\begin{equation}\label{eq12A}
M_p(f_{m_k}(\Gamma_{m_k}^1))=M_p(\Gamma_{m_k})\,.
\end{equation}
Однако, соотношения (\ref{eq11A}) и (\ref{eq12A}) в совокупности
противоречат (\ref{eq13A}). Полученное противоречие указывает на то,
что исходное предположение (\ref{eq6***A}) было неверным, и, значит,
семейство отображений $\frak{F}_{Q, \delta, K, n}(D, D^{\,\prime})$
равностепенно непрерывно в каждой точке $x_0\in \partial D.$
\end{proof}$\Box$

\medskip
{\bf 4. Доказательство основных результатов.} Утверждение теорем
\ref{th3} и \ref{th4} непосредственно вытекает из доказанных выше
лемм \ref{lem1} и \ref{lem1A} и \cite[лемма~3.1 и детали
доказательства теоремы~4.2]{GSS} (см. также
\cite[лемма~8]{Sev$_1$}).~$\Box$

{\bf 5. Несколько замечаний о точности условий.} Ограничимся для
простоты случаем $p=n.$ Прежде всего заметим, что в лемме \ref{lem1}
и теореме \ref{th3} нельзя, вообще говоря, отказаться от условия
наличия такого континуума $K,$ что $h(f^{-1}(K),
\partial D)\geqslant \delta>0,$ как
показывает простой пример семейства отображений $f(z)=z^n,$ $D=B(0,
1)\subset {\Bbb C}.$ Здесь указанное семейство отображений является
равностепенно непрерывным в $D,$ но не является равностепенно
непрерывным в $\overline{D},$ поскольку оно не является нормальным в
этой замкнутой области. Несколько сложнее построить пример семейства
$\{\frak{F}\}$ кольцевых $Q$-отображений, являющихся открытыми,
дискретными, удовлетворяющими условию $C(f,
\partial D)\subset \partial D^{\,\prime}$ для каждого
$f\in\{\frak{F}\},$ удовлетворяющих условию $h(f^{-1}(K),
\partial D)\geqslant \delta>0$ для некоторого континуума $K$ и числа $\delta>0,$
но при этом не являющегося равностепенно
непрерывным. Для этого зафиксируем числа $p\ge 1$ и $\alpha\in (0,
n/p(n-1)).$ Можно считать, что $\alpha<1$ в силу произвольности
выбора $p.$ Зададим последовательность гомеоморфизмов $g_m: {\Bbb
B}^n\setminus\{0\}\rightarrow {\Bbb R}^n$ следующим образом:
$$ g_m(x)\,=\,\left
\{\begin{array}{rr} \frac{1+|x|^{\alpha}}{|x|}\cdot x\,, & 1/m\le|x|\le 1, \\
\frac{1+(1/m)^{\alpha}}{(1/m)}\cdot x\,, & 0<|x|< 1/m \ .
\end{array}\right.
$$
Покажем, что последовательность $g_m$ удовлетворяет условиям
(\ref{eq3*!!})--(\ref{eq28*}) при некотором $Q(x).$ Заметим, что
каждое отображение $g_m$ переводит проколотый шар $D={\Bbb
B}^n\setminus \{0\}$ в кольцо $D^{\,\prime}=B(0,2)\setminus\{0\},$
которое, как известно, имеет сильно достижимую границу, что точка
$x_0=0$ является устранимой особенностью каждого $g_m,$ $m\in {\Bbb
N},$ причём $\lim\limits_{x\rightarrow 0}g_m(x)=0,$ и что
последовательность $g_m$ постоянна при $|x|\ge 1/m,$ а именно,
$g_m(x)\equiv g(x)$ при всех $x:\ \frac{1}{m}<|x|< 1,$
$m=1,2\ldots\,,$ где $g(x)=\frac{1+|x|^{\alpha}}{|x|}\cdot x.$
Заметим, что $g_m\in ACL({\Bbb B}^n).$ Действительно, отображения
$g_m^{(1)}(x)=\frac{1+(1/m)^{\alpha}}{(1/m)}\cdot x,$
$m=1,2,\ldots,$ являются отображениями класса $C^1,$ скажем, в шаре
$B(0, 1/m+\varepsilon)$ при малых $\varepsilon>0,$ а отображения
$g_m^{(2)}(x)=\frac{1+|x|^{\alpha}}{|x|}\cdot x$ ~--- отображениями
класса $C^1,$ скажем, в кольце $A(0, 1/m-\varepsilon, 1)=\left\{x\in
{\Bbb R}^n: 1/m-\varepsilon<|x|<1\right\}$ при малых
$\varepsilon>0.$ Отсюда вытекает, что гомеоморфизмы $g_m$ являются
липшицевыми в ${\Bbb B}^n$ и, значит, $g_m\in ACL({\Bbb B}^n),$ см.,
напр., \cite[разд.~5 на с.~12]{Va}. Далее, в каждой регулярной точке
$x\in D$ отображения $f:D\rightarrow {\Bbb R}^n$ рассмотрим {\it
внутреннюю дилатацию отображения $f$ в точке $x,$} определённую
равенством
$$K_I(x,f)\quad=\quad\frac{|J(x,f)|}{{l\left(f^{\,\prime}(x)\right)}^n}\,,$$
где $f^{\,\prime}(x)$ ~--- матрица Якоби отображения $f$ в точке
$x,$ $J(x, f)={\rm det\,}f^{\,\prime}(x)$ и
$l\left(f^{\,\prime}(x)\right)\,=\quad\min\limits_{h\in {\Bbb R}^n
\backslash \{0\}} \frac {|f^{\,\prime}(x)h|}{|h|}.$ Вычисляя
$K_I(x,f)$ для $f:=g_m,$  можно показать, что
$$ K_I(x, g_m)\,=\,\left
\{\begin{array}{rr} \left(\frac{1+|x|^{\,\alpha}}{\alpha
|x|^{\,\alpha}}\right)^{n-1}\,, & 1/m\le|x|\le 1, \\
1\,,\qquad & 0<|x|< 1/m \ ,
\end{array}\right.
$$
см.  \cite[предложение 6.3 гл. VI]{MRSY}. Заметим, что при каждом
фиксированном $m\in {\Bbb N},$ $K_I(x, g_m)\le c_m$ при некоторой
постоянной $c_m\ge 1.$ Значит, $g_m\in W_{loc}^{1, n}({\Bbb B}^n)$ и
$g_m^{-1}\in W_{loc}^{1, n}(B(0, 2)),$ поскольку условие $K_I(x,
g_m)\le c_m$ влечёт, что $g_m$ и $g_m^{-1}$ квазиконформны, см.,
напр., \cite[следствие 13.3 и теорема 34.6]{Va}. Тогда по
\cite[теорема~6.1 гл.~VI]{MRSY} гомеоморфизмы $g_m$ удовлетворяют в
области $D={\Bbb B}^n\setminus \{0\}$ неравенству вида
(\ref{eq3*!!}) при $Q=Q_m(x)= K_I(x, g_m).$
Более того, последовательность $g_m$ удовлетворяет неравенству вида
(\ref{eq3*!!}) при $Q=\left(\frac{1+|x|^{\,\alpha}}{\alpha
|x|^{\,\alpha}}\right)^{n-1}.$ 
Поскольку $\alpha p(n-1)<n,$ имеем, что $Q\in L^p({\Bbb B}^n).$ С
другой стороны, легко видеть, что
\begin{equation}\label{eq2!!!!!}
\lim\limits_{x\rightarrow 0}\ |g(x)| = 1\,,
\end{equation}
и $g$ отображает проколотый шар ${\Bbb B}^n\setminus\{ 0\}$ на
кольцо $1<|y|< 2.$ Тогда, в виду (\ref{eq2!!!!!}), мы получаем, что
$$|g_m(x)|=|g(x)|\ge 1\qquad\qquad\forall\quad x:|x|\ge 1/m,\quad m=1,2,\ldots\,,$$
т.е., семейство $\{g_m\}_{m=1}^{\infty}$ не является равностепенно
непрерывным в нуле. $\Box$

\medskip
КОНТАКТНАЯ ИНФОРМАЦИЯ

\medskip
\noindent{{\bf Евгений Александрович Севостьянов} \\
Житомирский государственный университет им.\ И.~Франко\\
кафедра математического анализа, ул. Большая Бердичевская, 40 \\
г.~Житомир, Украина, 10 008 \\ тел. +38 066 959 50 34 (моб.),
e-mail: esevostyanov2009@mail.ru}


\begin{thebibliography}{99}
{\small

\bibitem{Sev$_1$} {\it Севостьянов Е.А.} О точках ветвления отображений с неограниченной
характеристикой квазиконформности // Сиб. матем. ж. -- 2010. --
\textbf{51}, № 5. -- С. 1129--1146.

\bibitem{MRSY} {\it Martio O., Ryazanov V., Srebro U. and Yakubov
E.} Moduli in Modern Mapping Theory. -- New York: Springer Science +
Business Media, LLC, 2009.

\bibitem{Sev$_2$} {\it Севостьянов Е.А.} О равностепенной непрерывности гомеоморфизмов
с неограниченной характеристикой // Математические труды. -- 2012.
-- \textbf{15}, № 1. -- С. 178--204.


\bibitem{NP}{\it N\"{a}kki R. and Palka B.} Uniform equicontinuity
of  quasiconformal mappings // Proc. Amer. Math. Soc. -- 1973. --
\textbf{37}, no. 2. -- P. 427--433.

\bibitem{Vu} {\it Vuorinen M.} Exceptional sets and boundary behavior of quasiregular
mappings in $n$--space // Ann. Acad. Sci. Fenn. Ser. A 1. Math.
Dissertationes. -- 1976. -- \textbf{11}. -- P. 1--44.

\bibitem{Zel} {\it Зелинский Ю.Б.} Некоторые критерии гомеоморфизма
при отображении областей евклидова пространства // Труды VIII летней
математической школы. -- Киев: Ин-т математики АН УССР, 1971. -- С.
194--211.

\bibitem{Pol} {\it Полецкий Е.А.} Метод модулей для негомеоморфных квазиконформных
отображений // Матем. сб. -- 1970. -- \textbf{83,} № 2. -- С.
261--272.

\bibitem{Ri} {\it Rickman S.} Quasiregular mappings. -- Results in Mathematic and
Related Areas (3), 26. Berlin: Springer-Verlag, 1993.


\bibitem{Ku}{\it Куратовский К.} Топология, Т.\,2. -- М.: Мир, 1969.

\bibitem{He} {\it Heinonen~J.} Lectures on Analysis on metric spaces.
-- New York: Springer Science+Business Media, 2001.

\bibitem{Sev$_6$} {\it Sevost'yanov E.A.} On open and discrete
mappings with a modulus condition // Ann. Acad. Sci. Fenn. -- 2016.
-- \textbf{41}. -- P. 41–-50.

\bibitem{AS} {\it Adamowicz~T. and Shanmugalingam~N.}
Non-conformal Loewner type estimates for modulus of curve families
// Ann. Acad. Sci. Fenn. Math. -- 2010. -- \textbf{35.} -- P.~609–-626.

\bibitem{Sev$_4$} {\it Севостьянов~Е.А.}  О граничном поведении открытых дискретных
отображений с неограниченной характеристикой // Укр. матем. ж. --
2012. -- \textbf{64,} № 6. -- С. 855--859.

\bibitem{Smol} {\it Смоловая Е.С.} Граничное поведение кольцевых
$Q$-го\-ме\-о\-мор\-физ\-мов в метрических пространствах  // Укр.
матем. ж. -- 2010. -- \textbf{62}, № 5. -- С. 682--689.


\bibitem{Sev$_3$} {\it Севостьянов~Е.А.} Теория модулей, ёмкостей и нормальные семейства
отображений, допускающих ветвление // Украинский матем. вестник. --
2007. -- \textbf{4,} № 4. --  С. 582--604.


\bibitem{Va} {\it V\"{a}is\"{a}l\"{a} J.} Lectures on $n$-Dimensional Quasiconformal
Mappings. --  Lecture Notes in Math. 229, Berlin etc.:
Springer--Verlag, 1971.

\bibitem{Sev$_5$} {\it Севостьянов~Е.А.} О некоторых свойствах обобщённых
квазиизометрий с неограниченной характеристикой // Укр. матем. ж. --
2011. -- \textbf{63}, № 3. -- С. 385--398.

\bibitem{GSS} {\it Golberg~A., Salimov~R. and Sevost'yanov~E.} Singularities of
discrete open mappings with controlled $p$-module // J. Anal. Math.
-- 2015. -- \textbf{127.} -- P. 303--328.

 }

\end{thebibliography}
\end{document}